\documentclass{svproc}
\usepackage{url}

\usepackage{amsfonts}
\usepackage{amsmath, amssymb} 
\usepackage{graphicx, color}
\usepackage{slashbox}
\usepackage{multirow}

\newcommand{\wrighta}[3]{\ensuremath{ W \left(\left. #1,  #2 \right|  #3 \right) }}

\newcommand{\hsym}[2]{  \ensuremath{ \; _{#1}F_{#2} }  }

\newcommand{\ai}{ \ensuremath{ \operatorname{Ai} } }

\newcommand{\llim}[3]{\ensuremath{ \lim\limits_{ #1 \rightarrow #2} #3 }}
\newcommand{\fclass}[2]{\ensuremath{  \mathbb{#1}^{\, #2} }}
\newcommand{\erf}[1]{\ensuremath{ \mathrm{erf} \left( #1 \right)   }}
\newcommand{\erfc}[1]{ \ensuremath{ \operatorname{erfc}\left( #1\right)  } }

\begin{document}
\mainmatter              
\title{Computation of the Wright function from its integral representation}
\titlerunning{Computation of the Wright function}  
%
\author{Dimiter Prodanov\inst{1,}\inst{2}
 }

\authorrunning{D. Prodanov} 
%
\tocauthor{Dimiter Prodanov}
\institute{EHS and NERF, Interuniversity Microelectronics Centre (IMEC), 3001 Leuven, Belgium,\\
\email{ dimiter.prodanov@imec.be} \\
\and
ITSDP, Institute of Information and
Communication Technologies (IICT), Bulgarian Academy of Sciences, 1431 Sofia, Bulgaria}

\maketitle              

\begin{abstract}
The Wright function arises in the theory of the fractional differential equations. 
It is a very general mathematical object having diverse connections with other special and elementary functions. 
The Wright function provides a unified  treatment of several classes of special functions, such as the Gaussian, Airy, Bessel, error functions,  etc. 
The manuscript presents a novel numerical technique for approximation of the Wright function using quadratures.
The algorithm is implemented as a standalone library using the double-exponential quadrature integration technique using the method of stationary phase. 	
Function plots for a  variety of parameter values are demonstrated. 

\end{abstract}
%

\section{Introduction}\label{sec:intro}
The Wright function, introduced by E.M. Wright in two seminal publications \cite{Wright1940,Wright1935}, is a special mathematical function originally defined by the infinite series:
\begin{equation}
W\left(a,b \middle| \ z\right) := \sum\limits_{k=0}^{\infty} \frac{z^k}{k! \, \Gamma (a k + b)}, \quad z, b \in \mathbb{C}, \quad a>-1
\end{equation}
where $\Gamma$ denotes the Euler's gamma function.
If the parameter \textit{a} is a positive real number, it is classified as the Wright functions of the first kind, and when $-1 < a < 0$ -- as the Wright functions of the second kind \cite{Mainardi2010a}. 

For rational values of $a$ the Wright function is reducible to a finite sum of hypergeoemtric functions, which are amenable to automated representation by computer algebra systems.
Many more formulas are tabulated in the recent work of \cite{Apelblat2021}.
The analytical properties of the Wright function will not be repeated here in view of space limitations. 

Until recently efficient numerical algorithms for its approximation were not available. 
The seminal studies of Luchko et al. treated only the case $\left|b\right| \le 1$\cite{Luchko2008,Luchko2010}.
There is also a more recent interest in numerical algorithms for its computation based on  Laplace transform inversion \cite{Aceto2022}. 
In view of its broad applications, methods of its computation can be of general interest.

\section{Applications}\label{sec:applications}
The Wright function arises in the theory of the space-time fractional diffusion equation (FDE) with the temporal Caputo derivative \cite{Gorenflo2000,Lipnevich2010}.
We recall that the Caputo's fractional derivative of order $\beta > 0$ is defined for $\beta \notin \mathbb{N}$ as the differ-integral
\begin{equation}
\mathcal{D}^\beta_t f (t) = \frac{1}{\Gamma (m - \beta )} \int_{0}^{t} \frac{f^{(m)} (u) du}{(t-u)^{\beta+1 - m}}
\end{equation}
Notably, the Green function of the time FDE is expressed as the function $ \wrighta {-a}{1-a}{z} $ of the similarity variable $z$, named after Mainardi (M-Wright function). 
The fractional differential equation in the Caputo sense with variable coefficients
\begin{equation}
\mathcal{D}^\beta_t t^q f^\prime (t) = \beta t^{q-1} f(t)
\end{equation}
admits for a solution $ f(t) = \wrighta{\beta}{q}{t^\beta} $ \cite{Garra2021}. 

The subordination formulas for the fractional diffusion-wave equations can also be expressed in terms of Wright functions.  
Recent surveys about its applications can be found in \cite{Mainardi2020} and \cite{Povstenko2021}.

The Wright function provides a unified  treatment of several classes of special functions.
Noteworthy are the relationships of the Wright to the Bessel functions:  
\begin{flalign}
	\wrighta {1}{\nu+1} {- \frac{z^2 } {4} } & =\left( \frac{z}{2}\right)^{-\nu} J_\nu(z),  \\ 	
	\wrighta {1}{\nu+1} { \frac{z^2 } {4} } & =\left( \frac{z}{2}\right)^{-\nu} I_\nu(z)
\end{flalign}

\section{The canonical complex integral representation}\label{sec:intrep}

The canonical complex integral representation of the Wright function is given by the line integral
\begin{equation}
\wrighta{a}{b}{z}  = \frac{1}{2 \pi i} \int_{Ha^{-}} \frac{e^{\xi + z \xi^{-a}}  }{\xi^{b}} d \xi, \quad z \in \fclass{C}{}
\end{equation}
along the Hankel contour, which surrounds the negative real semi-axis. 
For integral values of $b$ and $a $ the path of integration closes around the origin \textit{O}
and can be used to extend the domain of the function into the negative integer parameters. 
\begin{figure}[pt]
	\centering
	\includegraphics[width=0.5\linewidth]{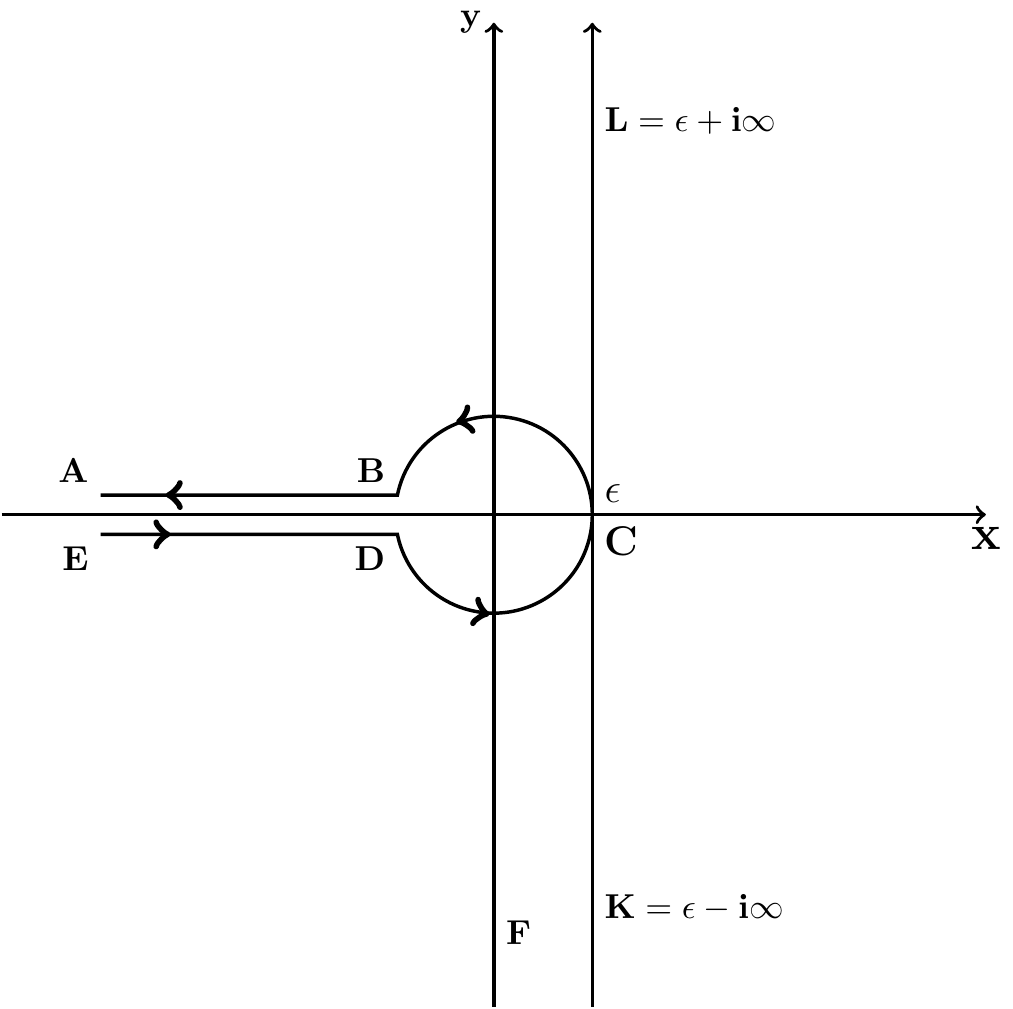}
	\caption{Hankel and Bromwich contours}
	\label{fig:hankel}
\end{figure}

\begin{definition}\label{def:kernel}
The integral kernel of the Wright function is defined as function of the complex variable $\xi$ as
\begin{equation}
ker(\xi):=\frac{e^{\xi + z \xi^{-a}}  }{\xi^{b}} 
\end{equation}
\end{definition}
\begin{definition}\label{def:holder}
The H\"older exponent function the origin is defined as
$
h(\xi):=\xi \frac{d}{d\xi}\log{\ker(\xi) } 	
$.
\end{definition}
\section{Derivation of the real integral representation}\label{sec:rep}
The results of this section are based on the method of Luchko \cite{Luchko2008}.
The starting point is a deformation of the Hankel contour into an  arc enclosing the pole at the origin and rays extending towards negative infinity as depicted in Fig. \ref{fig:hankel}.
The Hankel integral contour can be split in three parts - the rays $BA$, $ED$ and the arch $BCD$ (Fig. \ref{fig:hankel}).
Accordingly, the value of the integral can be evaluated as the sum:
\begin{multline}
\wrighta{a}{b}{z}= \underbrace{I_{AB}+I_{DE}}_{I_r} + \underbrace{I_{BCD}}_{P_\epsilon}
 = \\ \int_{AB} ker(\xi) d \xi + \int_{BCD} ker(\xi) d \xi + \int_{DE} ker(\xi) d \xi
\end{multline}
The branch cut will be taken along the negative real axis.
Then for the two rays $AB$ and $ED$ we obtain
\begin{equation}
	ker_{AB}= \frac{{{e}^{\frac{{{e}^{-i a \delta}} z}{{{r}^{a}}}+{{e}^{i \delta}} r-i b \delta}}}{{{r}^{b}}}, \quad
	ker_{ED}=\frac{{{e}^{\frac{{{e}^{i a \delta}} z}{{{r}^{a}}}+{{e}^{-i \delta}} r+i b \delta}}}{{{r}^{b}}}
\end{equation}
Therefore, 
\begin{equation}
	ker_{AB}-ker_{ED} = \frac{2 i {{e}^{\cos{(\delta)} r+\frac{\cos{\left( a \delta\right) } z}{{{r}^{a}}}}} \sin{\left( \frac{\sin{\left( a \delta\right) } z}{{{r}^{a}}}-\sin{(\delta)} r+b \delta\right) }}{{{r}^{b}}} 
\end{equation}
Taking the limit $\delta \to \pi $ results in the radial contribution
\begin{equation}
	I_r(\epsilon) = \frac{1}{ \pi} \int_{\epsilon}^{\infty} \frac{{{e}^{\frac{\cos{\left( \pi  a\right) } z}{{{r}^{a}}}-r}} \sin{\left( \frac{\sin{\left( \pi  a\right) } z}{{{r}^{a}}}+\pi  b\right) }}{{{r}^{b}}} dr
\end{equation}
For the arcuatic component we obtain
\begin{multline}
	ker_{BCD} =\frac{i {{e}^{\frac{{{e}^{-i a \varphi}} z}{{{\epsilon}^{a}}}+\epsilon\, {{e}^{i \varphi}}-i b \varphi}}}{{{\epsilon}^{b}}} e^{i \varphi} \epsilon\\
	 = \frac{i \epsilon {{e}^{\epsilon \cos{(\varphi)}+\frac{\cos{\left( a \varphi\right) } z}{{{\epsilon}^{a}}}}} \cos{\left( \frac{-\sin{\left( a \varphi\right) } z}{{{\epsilon}^{a}}}+\epsilon \sin{(\varphi)}-b \varphi +\varphi \right) }}{{{\epsilon}^{b}}} \\
	 -\frac{\epsilon {{e}^{\epsilon \cos{(\varphi)}+\frac{\cos{\left( a \varphi\right) } z}{{{\epsilon}^{a}}}}} \sin{\left( \frac{-\sin{\left( a \varphi\right) } z}{{{\epsilon}^{a}}}+\epsilon \sin{(\varphi)}-b \varphi + \varphi\right) }}{{{\epsilon}^{b}}} 
\end{multline}	
Therefore, taking  the limit  $\delta \to \pi $  as before we obtain the circular integral
\begin{equation}
	P_\epsilon = \frac{\epsilon^{1-b}}{2 \pi} \int_{-\pi}^{\pi} e^{\epsilon \cos{\varphi}+\cos{(a \varphi)} z/\epsilon^a} 
	\cos{(-( \sin(a\varphi) z/\epsilon^a+\epsilon \sin{\varphi}+ (1- b)\varphi )} d \varphi
\end{equation}
What is left is to determine the optimal value of the radial parameter $\epsilon$. At this point, the present approach departs from the treatment of \cite{Luchko2008}.

First, we can observe that if $a<0$ and  $b<1$ for $ { \epsilon\rightarrow 0 \Longrightarrow P_\epsilon \rightarrow 0}$ since the integral is bounded.
Furthermore, for $b=1$ and $a<0$
\begin{equation}
\llim{\epsilon}{0}{ P_\epsilon} = \frac{ 1}{2 \pi}\int_{-\pi}^{\pi}  \llim{\epsilon}{0}{} e^{\epsilon \cos{\varphi}+\cos{(a \varphi)} z/\epsilon^a} 
\cos{( -\sin(a\varphi) z/\epsilon^a+\epsilon \sin{\varphi} )} d \varphi= 1
\end{equation}
by Azrel\'a's theorem, which represents a useful special case. 
However, a more systematic way of determining $\epsilon$ will be the following. 

The H\"older exponent $\beta$ at the origin is determined by the limit of the  H\"older function
as $\beta =  \llim{\xi}{0}{h(\xi) }$.
For the Wright kernel one can directly  compute  
\begin{equation} 	 
h(\xi)=-\frac{a z}{{{\xi}^{a}}}+\xi-b 
\end{equation}
Therefore, for $a<0$ and  $ \beta= h(\xi)\Big|_{\xi=0} = -b $. This corresponds with the fact that   $\wrighta{a}{b}{0}=1/\Gamma(b)$.
Notably, if $b<1$ the singularity at the origin is integrable and $\epsilon=0$.

Furthermore, it can be concluded that for integral $a$ and $b$ the contour can be closed and the Cauchy theorem can be applied.
Let's denote $n=-a$ and $b=m$ for $ m, n \in \fclass{N}{}$; then
\begin{equation}
\wrighta{-n}{m}{z} = \operatorname{Res} \left( ker (\xi), \xi=0\right) = \frac{1}{\Gamma(m)} \left( \frac{d}{d \xi}\right)^{m-1}  e^{\xi + z \xi^{n}} \Bigg|_{\xi=0}
\end{equation}
where the differentiation takes precedence over the limit. 
Therefore, we can conclude that \wrighta{-n}{m}{z} is a polynomial in $z$.
This is a novel result, which was not anticipated by Wright and Mainardi, and allows for the extension of the domain of the parameters of the function. 
If $b=1$ then $ \operatorname{Res}  \left( ker (\xi), \xi=0\right) =1$.

In contrast, if  $a>0$ an essential singularity appears at the origin as the H\"older exponent diverges.
For that latter case we look for the curve where the phase of the kernel is stationary to first order to minimize oscillations. 
This is given by the computation
\begin{equation}
	\frac{d}{d\xi} \left( z/\xi^a+\xi \right)\Big|_{\xi=\xi_0} =\left(  -\frac{a z}{{{\xi}^{a+1}}}+1\right)  \Big|_{\xi=\xi_0} =0  \Longrightarrow \xi_0= \sqrt[a+1]{a z}
\end{equation}
Therefore, we can take
$   \epsilon=  |\xi_0| = \sqrt[a+1]{ |a z|} $.
\begin{remark}
	In particular cases (i.e. when $a>1$) it can be advantageous to change the integration variable as $r^a \mapsto u$
	\begin{equation}
		I_{u} (\epsilon)=\frac{1}{ \pi a } \int_{\epsilon}^{\infty} {u^{\frac{1-b}{a}-1}}{\sin{(\sin{(\pi a)} }z/u+\pi b)}  e^{\cos(\pi a) z/u-u^{1/a}} du 
	\end{equation}
\end{remark}
In summary, the following theorem can be stated.
\begin{theorem}
The Wright function can be computed according to  the table 

	\centering
	\begin{tabular}{|c|c|c|c|}
		\hline
		\backslashbox{b}{a}
		& $a<0$   & $a=0$   &   $a>0$\\
		\hline
		$b\leq1$	     &  $I_{r}(0) $   & 	\multirow{4}{*}{  $I_0 $}  &  \multirow{4}{*}{  $I_{u} (\epsilon^a)+ P (\epsilon)$} \\
		\cline{1-2}
		$b=1$ 		 & $I_{r}(0)+1$ &     &     \\
		\cline{1-2}
		$b>1$		 & $I_{r}(\epsilon) + P (\epsilon)$   &   &   \\
		\hline
	\end{tabular}

where $\epsilon = \operatorname{max} \left( \sqrt[a+1]{|a z|}, 1\right)$  and
\begin{flalign}
	I_0 &= {e^z}/{\Gamma(b)} \\
	I_r (\epsilon) & = \frac{1}{ \pi} \int\limits_{\epsilon}^{\infty} \frac{ \sin{\left( \frac{\sin{\left( \pi  a\right) } z}{r^{a}}+\pi  b\right) }}{{{r}^{b}}} {{e}^{\frac{\cos{\left( \pi  a\right) } z}{{{r}^{a}}}-r}}  dr   \\
	P (\epsilon) & = \frac{\epsilon^{1-b}}{2 \pi} \int\limits_{-\pi}^{\pi} e^{\epsilon \cos{\varphi}+\cos{(a \varphi)} z/\epsilon^a} 
	\cos{(\epsilon \sin{\varphi} -z\sin(a\varphi) /\epsilon^a+ (1- b)\varphi)} d \varphi,  
\end{flalign}
\end{theorem}

\subsection{The Double-Exponential (DE) quadrature method}\label{sec:de}
The DE quadrature integration technique of Takahasi and Mori \cite{Takahasi1973,Mori1985} can be summarized as follows.
The integral on the interval $A$ 
\begin{equation}
I= \int_{A} f(x) dx = \int_{A} f(x) dx= 
\int_{-\infty}^{\infty} f \circ \phi \ (t) \phi^{\prime} (t) dt 
\end{equation}
is transformed to an integral on the entire real line
and the function $\phi$ guarantees double exponential convergence to  \textit{I} as the limit of the Riemannian sum
\begin{equation}
	I= 	\llim{N}{\infty}{ } S_N =
	\llim{N}{\infty}{ }  h \sum\limits_{k=-N}^{N} \underbrace{\phi^{\prime}(k h)}_{w_k} f(\underbrace{\phi(k h)}_{x_k})  
\end{equation}
Therefore, the value of the integral can be approximated by the truncated sum $S_N \approx I$, where the weights are given by $w_k$ and the abscissas (i.e. evaluation points) by $x_k$, while $h$ is an adjustable parameter. 
For the case of a semi-finite interval the DE formula employs the transformation
\begin{equation}
 x= \phi(t)= \exp{\left( \frac{\pi}{2}\sinh{t}\right) } , \quad dx= \phi^\prime (t) dt= \frac{\pi }{2} {{e}^{\frac{\pi  \sinh{t}}{2}}} \cosh{t} \; dt
\end{equation}
For compact intervals the DE formula uses  the transformation
 \begin{equation}
	x=\phi(t)= \tanh{\left(\frac{\pi}{2}\sinh{t}\right)  }, \quad dx=\phi^\prime (t) dt= \frac{\pi }{2}   \operatorname{sech}{\left( \frac{\pi  \sinh{t}}{2} \right)}^2 \cosh{t}   dt
\end{equation}
The absolute error of the method is of the order $\exp{(–c\, N/log(N))}$. 

\section{Numerical Results}\label{sec:results}

Proposed algorithm is implemented as a standalone library using the DE quadrature integration technique \cite{Mori2001}  in the Java programming language and can be downloaded from \url{https://github.com/dprodanov/dspquad}.
A reference  implementation in the computer algebra system Maxima was also developed both for the QUADPACK \cite{piessens1983} and DE libraries.

Numerical experiments were performed on a 64-bit Microsoft Windows 10 Enterprise machine with configuration -- 
Intel\textsuperscript{\textregistered } Core\textsuperscript{TM} i5-8350U CPU @ 1.70GHz, 1.90 GHz and 16GB RAM.
Computation times are presented in Table \ref{tab:comptimes}, where $w_1=2.404825557695773$ is the first (approximate) zero of $J_0(x)$.
\begin{table}[h!]
	\centering
	\begin{tabular}{l|c|r|r}
		\hline
		Wright function & formula & x & time [$\mu$s] \\
		\hline
		\wrighta{-1/2}{1}{x}     & \erf{x/2}+1        & 2   & 204.400 \\
		\wrighta{-1/2}{1/2}{x}   & $M_{1/2}(-x)$      & 1/2 & 165.301 \\
		\wrighta{-1/2}{-1/2}{x}  & $G_2(x)$           & 3/2 & 199.200 \\
		\wrighta{-1/3}{2/3}{x}   & $ M_{1/3} (-x)$    & 1/2 & 187.500 \\
		\wrighta{1}{3/2}{x}      & $\ {\sin{\left( 2 \sqrt{x}\right) }}/{\sqrt{\pi x}}$ \ &  $\pi^2$ & 227.500 \\
		\wrighta{1}{1}{x}        &  $J_0(2 \sqrt{x})$ &\ $w_1^2/4$ & 250.301 \\
		\hline
	\end{tabular}
	\caption{Computation times}\label{tab:comptimes}
\end{table}
Plots were produced by Maxima 5.46.0 running Steel Bank Common Lisp (SBCL) v. 2.2.2.  
Presented results demonstrate that computation times are similar and do not depend much on the parameter values. 

\begin{figure}[h]
	\centering
	\begin{tabular}{ll}
		A & B \\
		\includegraphics[width=.5\textwidth]{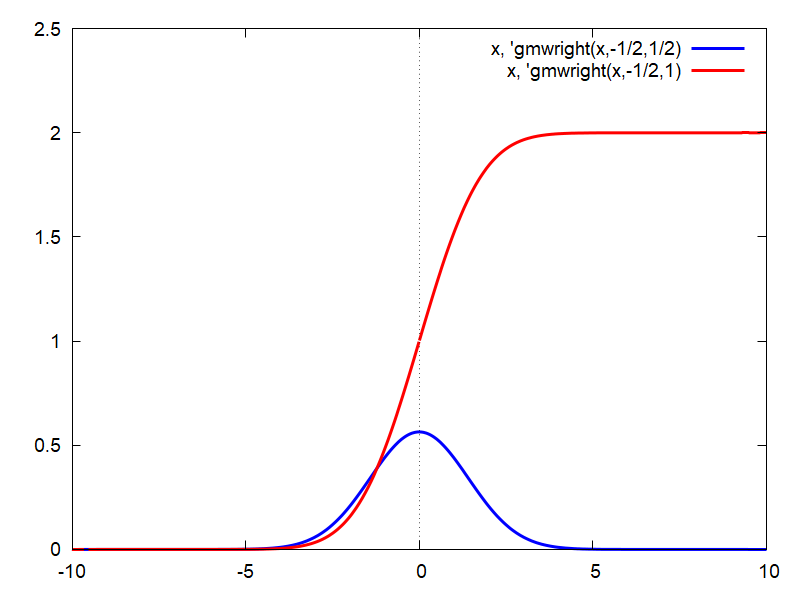} &
		\includegraphics[width=.5\textwidth]{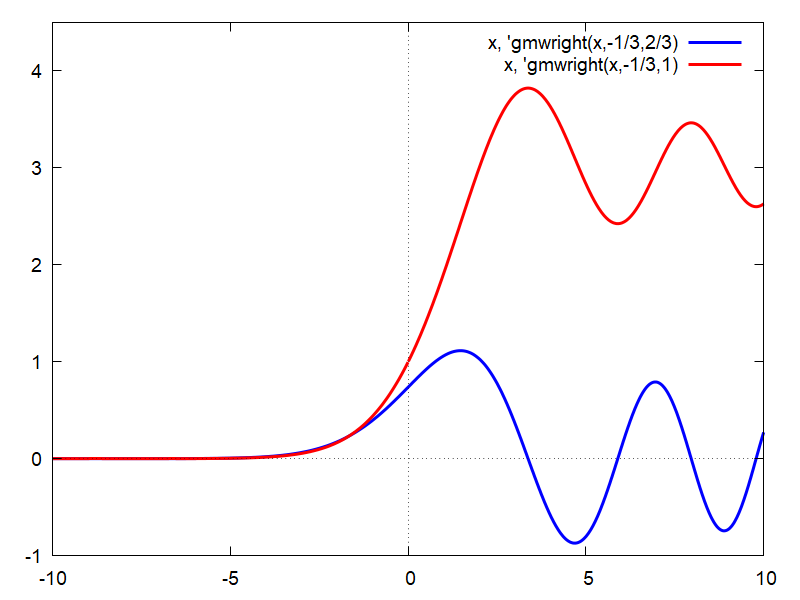} \\
	\end{tabular}
	\caption{Plots of the Wright function of the second kind, $a<0$:  
		A -- Gaussian and $\operatorname{erf}(x)$ functions; B -- $\operatorname{Ai}(x)$ and its integral }
	\label{fig:Mf}
\end{figure}
\begin{figure}[h]
	\centering
	\begin{tabular}{ll}
		A & B \\
		\includegraphics[width=.5\textwidth]{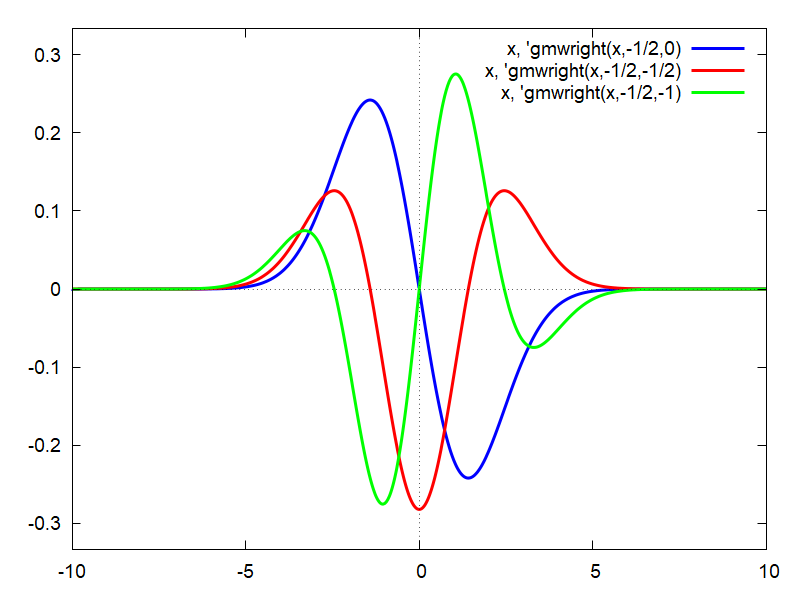} &
		\includegraphics[width=.5\textwidth]{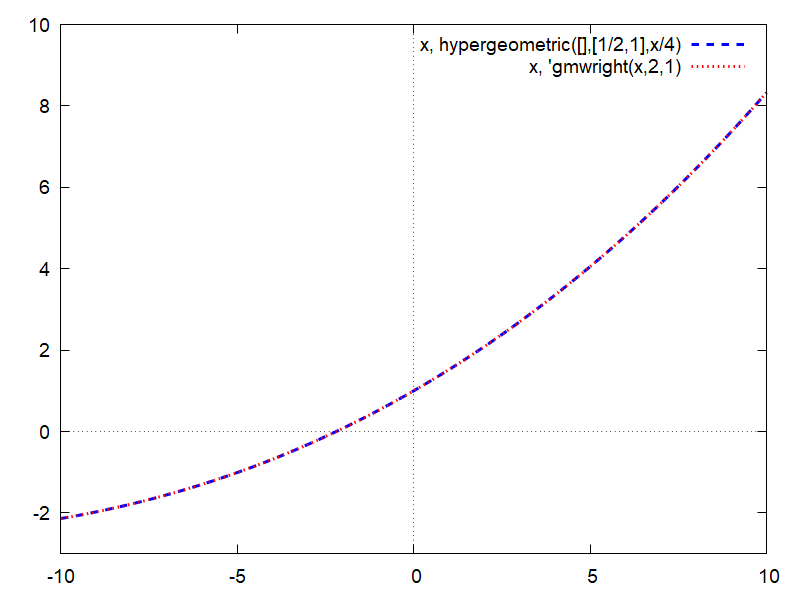} \\
	\end{tabular}
	\caption{ A -- Gaussian derivatives computed as Wright functions -- $G_1(x) \div G_3 (x)$;
		B -- hypergeometric function computed as a Wright function and overlaid on $\hsym{0}{2}(-; 1/2, 1| x/4)$ }
	\label{fig:gauss}
\end{figure}
\begin{figure}[h]
	\centering
	\centering
	\begin{tabular}{ll}
		A & B \\
		\includegraphics[width=.5\textwidth]{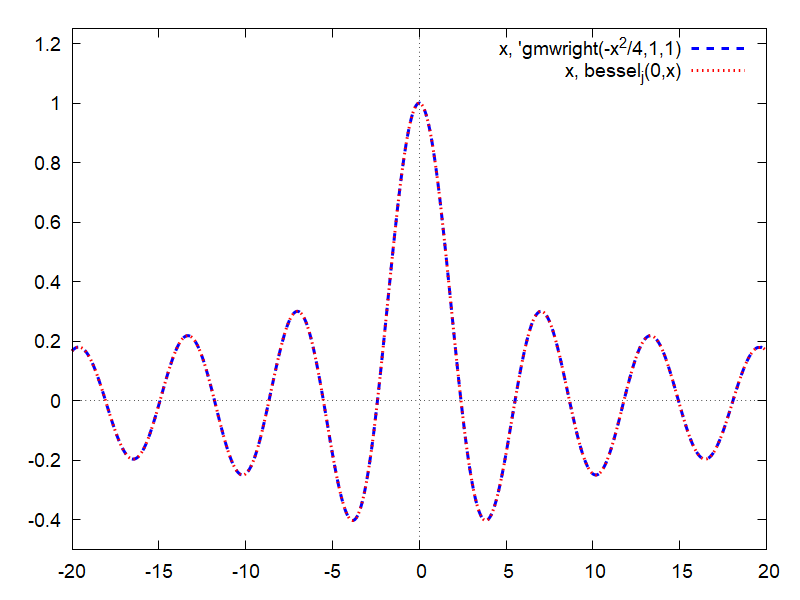} &
		\includegraphics[width=.5\textwidth]{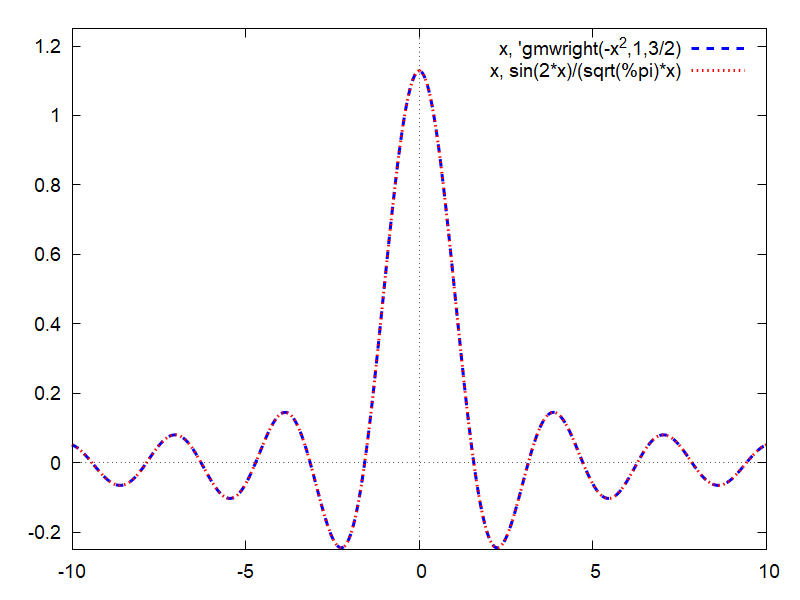} \\
	\end{tabular}
	\caption{Wright function plots for $a=1$: A -- $J_0 (x)$ overlaid on \wrighta{1}{1}{-x^2/4}, B -- $J_{1/2} (2 x)/\sqrt{|x|}$ overlaid on \wrighta{1}{3/2}{-x^2}  }
	\label{fig:bessel}
\end{figure}

\subsection{Computation of the Mainardi's Wright function and its integral}\label{sec:mwright}

Mainardi introduced a specialization of the Wright function,  called also M-Wright function, \cite{Mainardi2001}:
\begin{equation}
M_{a}( z) := \wrighta{- a}{1- a} {-z}
\end{equation}
with particular cases 
$
M_0 (z)=  e^{-z}$,   $M_{1/2} (z) =\frac{1}{\sqrt{\pi}} e^{-z^2 /4}$,   $M_{1/3} (z) =  \sqrt[3]{3^2} \ai \left( z/ \sqrt[3]{3}\right)
$.
For a rational $a=1/n$ the M-Wright function satisfies a differential equation of order $n-1$
\begin{equation}
\frac{d^{n-1}}{dz^{n-1}} M_{1/n}( z) + \frac{(-1)^n}{n} z M_{1/n}( z)=0
\end{equation}
equipped with the appropriate initial conditions. 
The integral of the M-Wright function is again a Wright function, but not a M-Wright function:
\begin{equation}
\int\limits_{-\infty}^{x} M_{a} (-u) du = \wrighta{- a}{1}{ x}
\end{equation}
A particular case of the Wright function is the error function:
\begin{equation}
	\wrighta{-\frac{1}{2}}{1}{z} = \erfc{-\frac{z}{2}}.
\end{equation}

Plots are presented in Fig. \ref{fig:Mf}A for $a=-1/2$, corresponding to the Gaussian and error functions, and in Fig. \ref{fig:Mf}B $a=-1/3$, corresponding to the Airy function and its integral.
%

\subsection{Computation of Gaussian derivatives and hypergeometric functions}\label{sec:gd}

The Wright function is also related to the derivatives of the Gaussian   and the Airy functions. 
The Gaussian derivatives can be represented as
\begin{equation}
G_n(x):=\left( \frac{d}{d x}\right)^n \frac{e^{-x^2/4}}{\sqrt{\pi}} =  \wrighta{- \frac{1}{2} }{\frac{1-n}{2}}{ x}
\end{equation}
Plots are presented in Fig. \ref{fig:gauss}A.
Recently, Apelblat and Gonzales-Santander \cite{Apelblat2021} established the equivalence
\begin{equation}
\wrighta{2}{1}{x} = \hsym{0}{2}\left( -, 1/2, 1; x/4\right) 
\end{equation}
A comparison plot is presented in Fig. \ref{fig:gauss}B.

\subsection{Computation of Bessel and trigonometric functions}\label{sec:bessel}
The algorithm shows good numerical stability for the computation of Bessel functions for intermediate values of the argument (Fig. \ref{fig:bessel}).
Bessel functions for half-integer arguments can be expressed in elementary trigonometric functions. For example, 
\begin{equation}
J_{1/2} ( x) = \frac{\sqrt{2} }{\sqrt{\pi |x|}} \sin{x}
\end{equation}
A comparison plot is presented in Fig. \ref{fig:bessel}B.


\section{Discussion}\label{sec:disc}
Computation of the Wright function presents an interesting applied mathematics problem because of its generality and its connections to other special functions, which so far have been computed in a variety of different ways. Here it is demonstrated that all of the above functions can be computed through the computation of the Wright function. 
The present contribution builds on the seminal works of Luchko et al. \cite{Luchko2008,Luchko2010} but removes some limitations found in their works (for example Th. 2.1. in \cite{Luchko2010}) and strives for full generality of the computational technique.
Unlike the algorithm of Luchko the present algorithm does not use series summation for small arguments, which saves on computation of $\Gamma$ functions. 
Unlike the method of Aceto and Durastante \cite{Aceto2022}, which uses a parabolic contour,  the present technique uses a combination of a  circular contour and a straight line integral.
Presented approach demonstrates very good computational times and will
allow for an easy extension to complex arguments.

\section*{Acknowledgment}
The present work is funded by the Horizon Europe project VIBraTE, grant agreement no. 101086815. 

\bibliographystyle{spmpsci_unsrt}
\bibliography{biodiffusion2}
 
\end{document}